\documentclass[12pt, a4paper]{article}
\usepackage[cp1251]{inputenc}
\usepackage[T2A]{fontenc}
\usepackage[russian]{babel}
\usepackage{amsmath}
\usepackage{float}
\usepackage{graphicx}
\usepackage{latexsym}
\usepackage{epsfig}
\usepackage{amsfonts}
\usepackage{amssymb}

\newtheorem{theorem}{Theorem}
\newtheorem{lemma}{Lemma}

\begin{document}

\title{~\\
[-6ex] Intermediately subcritical branching process in random environment:
the initial stage of the evolution\thanks{
This work is supported by the Russian Sciences Foundation under grant
19-11-00111 }}
\author{\textsc{Dyakonova E.E.
}\thanks{Steklov Mathematical Institute of Russian Academy of Sciences, 8 Gubkin St.,
Moscow, 119991, Russia; e-mail elena@mi.ras.ru} }
\date{}
\maketitle

\begin{abstract}
We consider branching process evolving in i.i.d. random environment. It is
assumed that the process is intermediately subcritical case. We investigate
the initial stage of the evolution of the process given its survival up to
the distant moment of time.
\end{abstract}



\section{Introduction and main result}

Consider a branching processes in random environment (BPRE). For the first
time BPRE have been introduced in \cite{at,sm} as a model for the
development of a population. Particles in this process reproduce
independently of each other according to some random reproduction laws which
can vary from one generation to the other. The complete and detailed
construction of the model of branching process in random environment is
given in the monograph \cite{VaKe17}. We give here a short description of
it. Denote by $\Delta $ the space of all probability measures on
$\mathbb{N}_{0}:=\{0,1,2,...,\}$. Equipped with a metric, it is a Polish space. Let $Q$
be a random variable taking values in $\Delta $. An infinite sequence $\Pi
=(Q_{1},Q_{2},\ldots )$ of i.i.d. copies of $Q$ is called a \textsl{random
environment } and $Q_{n}$ is the (random) offspring distribution of a
particles in generation $n-1$. Let $Z_{n}$ be the number of particles in
generation $n$.

A sequence of non-negative integervallued random variables
$Z_{0},Z_{1},\ldots $ is called a \emph{branching process in the random
environment }$\Pi $, if $Z_{0}$ is independent of $\Pi $\ and, given $\Pi $
the process $Z=(Z_{0},Z_{1},\ldots )$ is a Markov process with the law
\begin{equation}
\mathcal{L}\big(Z_{n}\;\big|\;Z_{n-1}=z,\,\Pi =(q_{1},q_{2},\ldots )\big)\
=\ q_{n}^{\ast z}  \label{transition}
\end{equation}%
for every $n\in \mathbb{N}=\{1,2,\ldots \}$, $z\in \mathbb{N}_{0}$ and
$q_{1},q_{2},\ldots \in \Delta $, where $q^{\ast z}$ is the $z$-fold \
convolution of the measure $q$.

We denote by $\mathbb{P}$ the corresponding probability measure on the
underlying probability space. For convenience we assume that $Z_{0}=1$ and
$\mathbb{P}(Q(0)=1)< 1$ (here and in what follows we use the notation
$Q(y),q(y)$ for $Q(\{y\}),q(\{y\})$).\medskip \newline
For a probability measure $q\in \Delta $ we define its mean \begin{equation*}
m(q):=\sum_{y=0}^{\infty }y\ q(y).
\end{equation*}%
Introduce the so-called \textsl{associated random walk }$S=(S_{n})_{n\geq
0}. $ This random walk has increments given by \begin{equation*}
\ X_{n}\ :=\ \log m(Q_{n}),\text{ \ }n\in \mathbb{N}_{0}.\
\end{equation*}These increments are i.i.d. probabilistic copies of $X\ :=\ \log m(Q).\ $
Assume that $S_{0}=0$ .

A BPRE is called \emph{critical} if $S$ is an oscillating random walk. A
BPRE is called \emph{subcritical} if its oscillating random walk drifts to
$-\infty $ $\mathbb{P}$--a.s. as $n\rightarrow \infty $. \ A BPRE is
\emph{weakly} subcritical if $\mathbb{E}[Xe^{X}]>0$, \emph{intermediately }
subcritical if $\mathbb{E}[Xe^{X}]=0,$ and \emph{strongly} subcritical if
$\mathbb{E}[Xe^{-X}]<0$.

Clearly  \begin{equation*}
\mathbb{E}[Z_{n}\,|\,\Pi \,]\ =\ \prod_{k=1}^{n}m(Q_{k})\ =\ \exp
(S_{n})\quad \mathbb{P}\text{--a.s.,}
\end{equation*}and as $n\rightarrow \infty $ \begin{equation*}
\mathbb{P}(Z_{n}>0\,|\,\Pi )\leq \min_{0\leq k\leq n}\mathbb{E}[Z_{k}|\Pi
]\leq \exp \big(\min_{0\leq k\leq n}S_{k}\big)\rightarrow 0\quad
\mathbb{P}\text{--a.s.}
\end{equation*}%
for critical and subcritical BPRE. We will assume that the following
restriction is valid.

\paragraph{Assumption A1.}

\emph{The process $Z$ is intermediately subcritical, i.e.} \begin{equation*}
\mathbb{E[}X]\ <\ 0,\;\mathbb{E}[Xe^{X}]\ =\ 0\ .
\end{equation*}

\noindent One of the essential instrument for the study of the properties of
PBRE is a change of measure. We also use this approach and introduce along
with the measure $\mathbb{P}$, one more measure $\mathbf{P,}$ setting for
every $n\in \mathbb{N}$ and every bounded measurable function $\varphi
:\Delta ^{n}\times \mathbb{N}_{0}^{n+1}\rightarrow \mathbb{R}$
\begin{equation}
\mathbf{E}[\varphi (Q_{1},\ldots ,Q_{n},Z_{0},\ldots ,Z_{n})]\ =\ \gamma
^{-n}\mathbb{E}\big[\varphi (Q_{1},\ldots ,Q_{n},Z_{0},\ldots
,Z_{n})e^{S_{n}-S_{0}}\big]\ ,  \label{iks1}
\end{equation}where \begin{equation*}
\gamma \ =\ \mathbb{E}[e^{X}]\ .
\end{equation*}%
(We include $S_{0}$ in this expression$,$ since later on we will also
consider cases where $S_{0}\neq 0$.) Note that under this change of measure
the relation $\mathbb{E}[Xe^{X}]=0$ gives \begin{equation}
\mathbf{E}[X]\ =\ 0\ .  \label{criticality}
\end{equation}%
Thus, the random walk $S$ is recurrent under $\mathbf{P}$. Note, that
(\ref{iks1}) yields \begin{equation*}
\mathbb{E}[Z_{n}]=\gamma ^{n}\ .
\end{equation*}

We also need the following assumption.

\paragraph{Assumption A2.}

\emph{The distribution of $X$ is absolutely continious and belongs with
respect to $\mathbf{P}$ to the domain of attraction of a stable law with
index $\alpha \in (1,2].$ } \newline

\noindent It follows from (\ref{criticality}) that there is an increasing
sequence of positive numbers \begin{equation*}
a_{n}\ =\ n^{1/\alpha }\ell _{n},
\end{equation*}%
where $\ell _{1},\ell _{2},\ldots $ is a sequence slowly varying at infinity
such that as $n\rightarrow \infty $ the sequence$\ \tfrac{1}{a_{n}}S_{n}$
weakly converges, i.e. \begin{equation*}
\mathbf{P}\big(\tfrac{1}{a_{n}}S_{n}\in dx\,\big)\ \rightarrow \ s(x)\,dx,
\end{equation*}%
where $s(x)$ is the density of the limiting stable law. Note that in the
case of finite variance $\sigma ^{2}=\mathbf{E}[X^{2}]<\infty $ we have
$\ell _{n}=\sigma $.

We also need the following moment restriction on the distibution of the
random variable \begin{equation*}
\zeta (a):=\ \frac{1}{m(Q)^{2}}\sum_{y=a}^{\infty }y^{2}Q(y)\ ,\quad a\in
\mathbb{N}\ .
\end{equation*}

\paragraph{Assumption A3.}

There exist \emph{\ $\varepsilon >0$ }and \emph{\ $a\in \mathbb{N},$ }such
that $\mathbf{E}[(\log ^{+}\zeta (a))^{\alpha +\varepsilon }]\ <\ \infty ,$
where $\log ^{+}x=\log (x\vee 1)$.

Let \begin{equation*}
\tau _{n}=\min \{k\leq n\mid S_{k}\leq S_{0},S_{1},\ldots ,S_{n}\}
\end{equation*}%
be the moment, when $S$ takes its minimum for the first time on the interval
$[0,n].$

Let $r_{n}\in \mathbb{N}_{0},n=1,2,...,$ and $\ r_{n}\rightarrow \infty
,r_{n}=o\left( n\right) ,n\rightarrow \infty .$ For brevity we will use the
notation\begin{equation*}
r=r_{n},\;\tau =\tau _{r}.
\end{equation*}

The main result of this paper is the following statement in which the symbol
$\overset{d}{\rightarrow }$ denotes weak convergence.

\begin{theorem}
\label{theomain} If Assumptions A1 -- A3 are valid and $r=r_{n}=o\left(
n\right) $ as $n\rightarrow \infty ,$ then

1) there is a random variable $\varkappa _{1}$ with values in $\mathbb{N}$
such that as $n\rightarrow \infty $ \begin{equation*}
\big(Z_{\tau _{r}}\mid Z_{n}>0\big)\overset{d}{\rightarrow }\varkappa _{1};
\end{equation*}%
2) there is a positive random variable $\varkappa _{2}$ such that as
$n\rightarrow \infty $ \begin{equation*}
\Big(\frac{Z_{r}}{e^{S_{r}-S_{\tau _{r}}}}\ \big|\
Z_{n}>0\Big)\overset{d}{\rightarrow }\varkappa _{2}.
\end{equation*}
\end{theorem}

Observe, that Theorem \ref{theomain} characterizes the behaviour of an
intermediatly subcritical branching process $Z$ in random environment at the
initial stage of the evolution of the process given its survival up to time
$n\rightarrow \infty .$ This result complements the results of
\cite{ABKVINTER}, which investigates the behaviour of a subcritical branching
process $Z$ in random environment at moments $r_{n}=t_{i}n,t_{i}\in \lbrack
0,1],i=1,...,k,$ under the condition of the survival of the process up to
$n\rightarrow \infty $ (see for this matter the remark after Lemma
\ref{theomain}). Note also that the distribution of the number of particles at
the initial period of the\ evolution for critical and weakly subcritical
BPRE given their survival up to a distant moment were investigated in
\cite{VaDy17} and \cite{VaDy19}.

\section{ Auxiliary results}

In this section we assume that the associated random walk $S$ can start \
from arbitrary value $S_{0}=x,$ and use the symbols $\mathbf{P}_{x}(\cdot )$
and $\mathbf{E}_{x}[\cdot ]$ to denote the corresponding probabilities and
expectations. Thus $\mathbf{P}=\mathbf{P}_{0}$.

For $n\geq 1$ set \begin{equation*}
L_{n}\ :=\ \min (S_{1},\ldots ,S_{n})\ ,\quad M_{n}\ :=\ \max (S_{1},\ldots
,S_{n}).
\end{equation*}

Introduce the renewal functions $u:\mathbb{R}\rightarrow \mathbb{R}$ and
$v:\mathbb{R}\rightarrow \mathbb{R},$ defined by \begin{align*}
u(x)\ & =\ 1+\sum_{k=1}^{\infty }\mathbf{P}(-S_{k}\leq x,M_{k}<0)\ ,\quad
x\geq 0\ ;\quad u(x)=0,\quad x<0\ , \\
v(x)\ & =\ 1+\sum_{k=1}^{\infty }\mathbf{P}(-S_{k}>x,L_{k}>0)\ ,\quad x<0\
;\quad v(x)=0,\quad x>0\ , \\
v(0)\ & =\mathbf{E}[v\left( X\right) ;X<0].
\end{align*}

It is known \cite{ABKVINTER}, that \begin{equation}
\begin{array}{rl}
\mathbf{E}[u(x+X);X+x\geq 0]\ =\ u(x)\ , & x\geq 0\ , \\
\mathbf{E}[v(x+X);X+x<0]\ =\ v(x)\ , & x\leq 0\ ,\end{array}
\label{harm}
\end{equation}%
that allwos to specify two new measures $\mathbf{P}^{+}$ and
$\mathbf{P}^{-}$. The construction procedure of these measures is standard and explained in
detail, for example, in \cite{agkv}, \cite{bedo}. We give a sketch of this
construction.

Consider the filtration $\mathcal{F}=(\mathcal{F}_{n})_{n\geq 0}$, where
$\mathcal{F}_{n}=\sigma (Q_{1},\ldots ,Q_{n},Z_{0},\ldots ,Z_{n})$. \ Thus,
$S $ is adapted to $\mathcal{F}$, and $X_{n+1}$ (as well as the measure
$Q_{n+1} $) is independent of $\mathcal{F}_{n}$ for all $n\geq 0$. For every
bounded, $\mathcal{F}_{n}$-measurable random variable $R_{n}$ set
\begin{align*}
\mathbf{E}_{x}^{+}[R_{n}]& :\ =\
\frac{1}{u(x)}\mathbf{E}_{x}[R_{n}u(S_{n});L_{n}\geq 0]\ ,\quad x\geq 0\ , \\
\mathbf{E}_{x}^{-}[R_{n}]\ & :=\
\frac{1}{v(x)}\mathbf{E}_{x}[R_{n}v(S_{n});M_{n}<0]\ ,\quad x\leq 0\ .
\end{align*}%
In other words, $\mathbf{P}_{x}^{+}$ and $\mathbf{P}_{x}^{-}$ correspond to
conditioning the random walk $S$, not to enter $(-\infty ,0)$ and $[0,\infty
),$ respectively.

Note that by duality (see \cite{Fel2}, ch. XII, \S 2)\ \begin{equation*}
\mathbf{P}(\tau _{n}=n)=\mathbf{P}(M_{n}<0)\ .
\end{equation*}We need the following results.

\begin{lemma}
\label{maxim} (\cite{ABKVINTER} ) Under Assumption A2 there are real numbers
\begin{equation*}
b_{n}=n^{1-\alpha ^{-1}}\ell _{n}^{\prime }\ ,\quad n\geq 1,
\end{equation*}%
where $(\ell _{n}^{\prime },n\geq 1)$ is a sequence slowly varying at
infinity, such that for every $x\geq 0$ as $n\rightarrow \infty $
\begin{equation*}
\mathbf{P}(M_{n}<x)\ \sim \ v(-x)b_{n}^{-1}.
\end{equation*}%
Moreover, there is a constant $c>0$ such that for all $x\geq 0$
\begin{equation*}
\mathbf{P}(M_{n}<x)=\mathbf{P}_{-x}(M_{n}<0)\leq cv(-x)b_{n}^{-1}\ .
\end{equation*}
\end{lemma}

\begin{lemma}
\label{limitEminus}(\cite{ABKVINTER}, Lemma 2.3) Assume A2 and let
$U_{1},U_{2},\ldots $ be a sequence of uniformly bounded random variables,
adapted to the filtration $\mathcal{F}$. If $U_{n}\rightarrow U_{\infty }$
$\mathbf{P}^{+}$-a.s. for some limiting random variable $U_{\infty }$, then
as $n\rightarrow \infty $ \begin{equation*}
\mathbf{E}[U_{n}\mid L_{n}\geq 0]\ \rightarrow \ \mathbf{E}^{+}[U_{\infty
}]\ .
\end{equation*}%
Similarly, if $U_{n}\rightarrow U_{\infty }$ $\mathbf{P}^{-}$-a.s., then as
$n\rightarrow \infty $ \begin{equation*}
\mathbf{E}[U_{n}\mid M_{n}<0]\ \rightarrow \ \mathbf{E}^{-}[U_{\infty }]\ .
\end{equation*}
\end{lemma}

\noindent Let $D[0,1]$ be the space of cadlag functions on $[0,1]$. If we
equippe space $D[0,1]$ with the Skorohod metric we get the Skorohod space
$D[0,1].$ \ For convinious we assume that $S_{nt}=S_{\lfloor nt\rfloor }$. It
follows from Assumption A2 (see, for example, \cite{bi}), that there exists
a Levy-process $L=(L_{t})_{0\leq t\leq 1}$ such that the process
$S^{n}:=(\tfrac{1}{a_{n}}S_{nt})_{0\leq t\leq 1}$ converges in disrtibution to $L$ in
the Skorohod space $D[0,1]$. \ Let $L^{-}=(L_{t}^{-})_{0\leq t\leq 1}$ be
the corresponding non-positive Levy meander. This is the process
$(L_{t})_{0\leq t\leq 1}$, conditioned on the event $\sup_{t\leq 1}L_{t}\leq
0 $ (see \cite{chau} and \cite{be}). Let $\overset{d}{\Longrightarrow }$
denote the weak convergence in the Skorohod space $D[0,1]$.

We need the following result.

\begin{lemma}
\label{funclimit2}(\cite{ABKVINTER}, Lemma 2.4) Under Assumptions A1 and A2
for every $x\geq 0$ as $n\rightarrow \infty $ \begin{equation*}
\big(\tfrac{1}{a_{n}}S^{n}\mid M_{n}<x\big)\overset{d}{\Longrightarrow }L^{-}
\end{equation*}in the Skorohod space $D[0,1]$.
\end{lemma}

\noindent Recall that the random walk $S$ satisfies, under the measure
$\mathbf{P,}$ Assumption A2. Let $S^{\left( 1\right) }$ and $S^{\left(
2\right) }$\ be two independent probability copies of the random walk
$S.$\thinspace\ For $S^{\left( 2\right) }$ we denote by $\tau _{n}^{\left(
2\right) }$and $M_{n}^{\left( 2\right) }$ analogies of the random variables
$\tau _{n}$ and $M_{n}.$

\begin{lemma}
\label{vmestoLem2.5} Under Assumptions A1 and A2 \begin{equation*}
\mathbf{P}\left( \min_{1\leq i\leq r_{n}}\left( S_{i}^{\left( 1\right)
}-S_{r_{n}}^{\left( 1\right) }\right) \leq S_{n}^{\left( 2\right) }\mid \tau
_{n}^{\left( 2\right) }=n\right) \rightarrow 0,\;n\rightarrow \infty .
\end{equation*}
\end{lemma}

\textbf{Proof. }Using duality and taking into account that $r=r_{n}=o\left(
n\right) ,$ we have \ \begin{eqnarray*}
&&\mathbf{P}\left( \frac{\min_{1\leq i\leq r}\left( S_{i}^{\left( 1\right)
}-S_{r}^{\left( 1\right) }\right) }{a_{n}}\leq \frac{S_{n}^{\left( 2\right)
}}{a_{n}};\tau _{n}^{\left( 2\right) }=n\right) /\mathbf{P}\left( \tau
_{n}^{\left( 2\right) }=n\right) \\
&=&\mathbf{P}\left( \frac{\min_{1\leq i\leq r}\left( S_{i}^{\left( 1\right)
}-S_{r}^{\left( 1\right) }\right) }{a_{n}}\leq \frac{S_{n}^{\left( 2\right)
}}{a_{n}};M_{n}^{\left( 2\right) }<0\right) /\mathbf{P}\left( \tau
_{n}^{\left( 2\right) }=n\right) \\
&=&\mathbf{P}\left( \frac{\min_{1\leq i\leq r}\left( S_{i}^{\left( 1\right)
}-S_{r}^{\left( 1\right) }\right) }{a_{n}}\leq \frac{S_{n}^{\left( 2\right)
}}{a_{n}}\mid M_{n}^{\left( 2\right) }<0\right) \frac{\mathbf{P}\left(
M_{n}^{\left( 2\right) }<0\right) }{\mathbf{P}\left( \tau _{n}^{\left(
2\right) }=n\right) } \\
&=&\mathbf{P}\left( \frac{\min_{1\leq i\leq r}\left( S_{i}^{\left( 1\right)
}-S_{r}^{\left( 1\right) }\right) }{a_{n}}\leq \frac{S_{n}^{\left( 2\right)
}}{a_{n}}\mid M_{n}^{\left( 2\right) }<0\right) \rightarrow 0,\;n\rightarrow
\infty ,
\end{eqnarray*}%
since $\min_{1\leq i\leq r}\left( S_{i}^{\left( 1\right) }-S_{r}^{\left(
1\right) }\right) /a_{n}\overset{d}{\rightarrow }0$\ as $n\rightarrow \infty
,r=o\left( n\right) ,$ and given $M_{n}^{\left( 2\right) }<0$ the
distribution of $S_{n}^{\left( 2\right) }/a_{n}$ converges weakly to an
absolutely continious distribution on $[0,\infty )$.

The lemma is proved.

Let $Q_{j}=Q_{1}$ for $j\leq 0$.

\begin{lemma}
\label{leQ} If Assumptions A1 and A2\ are valid, then for any $m\geq 0,k\geq
1$ given the event $\tau _{n-m}=n-m$ for $n\rightarrow \infty $ distribution
of \begin{equation*}
\Big(\big(Q_{\tau _{r}+1},\ldots ,Q_{\tau _{r}+k}\big),\big(Q_{\tau
_{r}},\ldots ,Q_{\tau _{r}-k+1}\big),\left( \frac{S_{\tau
_{r}}}{a_{r}},\frac{S_{n-m}}{a_{n}}\right) \Big)
\end{equation*}%
converges weakly to a probability measure $\mu _{k}^{+}\otimes \mu
_{k}^{-}\otimes \mu $, where $\mu _{k}^{+}$ and $\mu _{k}^{-}$ are the
distributions of $(Q_{1},\ldots ,Q_{k})$ under the probability measures
$\mathbf{P}^{+}$ and $\mathbf{P}^{-},$ respectively, and $\mu $ \ is
probability measure on $\mathbb{R}^{2}$.
\end{lemma}

\textbf{Proof.} For $l\geq 0$ set \begin{equation*}
Q^{+}(l):=(Q_{l+1},...,Q_{l+k})\ ,\quad Q^{-}(l):=(Q_{l},..,Q_{l-k+1})\ .
\end{equation*}

Let $\phi _{1},\phi _{2}:\Delta ^{k}\rightarrow \mathbb{R}$ be bounded
functions and $\phi _{3},\phi _{4}:\mathbb{R}\rightarrow \mathbb{R}$ be
bounded continuous functions. We have \begin{align}
& \mathbf{E}\big[\phi _{1}(Q^{-}(\tau _{r}))\phi _{2}(Q^{+}(\tau _{r}))\phi
_{3}(\tfrac{S_{\tau _{r}}}{a_{r}})\phi _{4}(\tfrac{S_{n-m}-S_{\tau
_{r}}}{a_{n}});\tau _{n-m}=n-m\big]  \notag \\
& =\sum_{l=0}^{r}\mathbf{E}\big[\phi _{1}(Q^{-}(l))\phi _{2}(Q^{+}(l))\phi
_{3}(\tfrac{S_{l}}{a_{r}})\phi _{4}(\tfrac{S_{n-m}-S_{l}}{a_{n}});\tau
_{r}=l,\tau _{n-m}=n-m\big]\ .  \label{Dec1}
\end{align}Introduce the random walk \begin{equation*}
S_{k}^{^{\prime }}=S_{k+l}-S_{l},k\geq 0.
\end{equation*}%
We will use prime to denote the random variables related to this random
walk. It is easy to see that\begin{align}
& \mathbf{E}[\phi _{1}(Q^{-}(l))\phi _{2}(Q^{+}(l))\phi
_{3}(\tfrac{S_{l}}{a_{r}})\phi _{4}(\tfrac{S_{n-m}-S_{l}}{a_{n}});\tau _{r}=l,\tau _{n-m}=n-m]
\notag \\
& =\mathbf{E}[\phi _{1}(Q^{-}(l))\phi _{3}(\tfrac{S_{l}}{a_{r}});\tau
_{l}=l]\mathbf{E}[\phi _{2}(Q^{\prime +}(0))\phi _{4}(\tfrac{S_{n-m-l}^{\prime
}}{a_{n}});L_{r-l}^{\prime }\geq 0,\tau _{n-m-l}^{\prime }=n-m-l]  \notag \\
& =\mathbf{E}[\phi _{1}(Q^{+}(0))\phi
_{3}(\tfrac{S_{l}}{a_{r}});M_{l}<0]\mathbf{E}[\phi _{2}(Q^{\prime +}(0))\phi _{4}(\tfrac{S_{n-m-l}^{\prime
}}{a_{n}});L_{r-l}^{\prime }\geq 0,\tau _{n-m-l}^{\prime }=n-m-l]\ .
\label{zero}
\end{align}

\bigskip Note that for $l>k$ \begin{align*}
& \frac{\mathbf{E}\big[\phi _{1}(Q^{+}(0))\phi
_{3}(\tfrac{S_{l}}{a_{r}});M_{l}<0\big]}{\mathbf{P}(M_{l}<0)} \\
& \quad =\mathbf{E}\big[\phi _{1}(Q^{+}(0))\mathbf{E}_{S_{k}}[\phi
_{3}(\tfrac{S_{l-k}}{a_{r}})\mid
M_{l-k}<0]\tfrac{\mathbf{P}_{S_{k}}(M_{l-k}<0)}{\mathbf{P}(M_{l}<0)};M_{k}<0\big]\ .
\end{align*}%
Therefore by Lemmas \ref{maxim}, \ref{funclimit2} and the dominated
convergence theorem, if $l=l_{n}\sim tr_{n}$ for some $0<t<1$, then
$a_{l_{n}}/a_{r_{n}}\sim t^{\frac{1}{\alpha }}\ $and \begin{align}
& \frac{\mathbf{E}\big[\phi _{1}(Q^{+}(0))\phi
_{3}(\tfrac{S_{l_{n}}}{a_{r}});M_{l_{n}}<0\big]}{\mathbf{P}(M_{l_{n}}<0)}  \notag \\
& \qquad \qquad \qquad \rightarrow \mathbf{E}\big[\phi
_{1}(Q^{+}(0))v(S_{k});M_{k}<0\big]\,\mathbf{E}[\phi _{3}(t^{\frac{1}{\alpha
}}L_{1}^{-})]  \notag \\
& \qquad \qquad \qquad =\mathbf{E}^{-}\big[\phi
_{1}(Q^{+}(0))\big]\,\mathbf{E}[\phi _{3}(t^{\frac{1}{\alpha }}L_{1}^{-})]\ .  \label{perv}
\end{align}

\bigskip Now consider the expectation \begin{equation*}
\mathbf{E}[\phi _{2}(Q^{\prime +}(0))\phi _{4}(\tfrac{S_{n-m-l}^{\prime
}}{a_{n}});L_{r-l}^{\prime }\geq 0,\tau _{n-m-l}^{\prime }=n-m-l]\ .
\end{equation*}%
Set $R=r-l,N=n-m-l.$ Inroduce the random walk $S_{i}^{\prime \prime
}:=S_{R+i}^{\prime }-S_{R}^{\prime },$ $i\geq 0.$ Let\begin{equation*}
A_{0}:=\left\{ S_{1}^{\prime }\geq 0,S_{2}^{\prime }\geq 0,...,S_{R}^{\prime
}\geq 0\right\} ,
\end{equation*}\begin{equation*}
A_{1}:=\left\{ S_{1}^{\prime \prime }\geq S_{N-R}^{\prime \prime
},S_{2}^{\prime \prime }\geq S_{N-R}^{\prime \prime },...,S_{N-R-1}^{\prime
\prime }\geq S_{N-R}^{\prime \prime }\right\} ,
\end{equation*}\begin{equation*}
B:=\left\{ S_{N-R}^{\prime \prime }<-S_{R}^{\prime }\right\}
,\;\overline{B}=\left\{ S_{N-R}^{\prime \prime }\geq -S_{R}^{\prime }\right\} .
\end{equation*}Clearly, the events $A_{0}$ and$A_{1}$ are independent and
\begin{equation*}
\left\{ L_{r-l}^{\prime }\geq 0,\tau _{n-m-l}^{\prime }=n-m-l\right\}
=A_{0}\cap A_{1}\cap B.
\end{equation*}Therefore, \begin{align}
& \mathbf{E}[\phi _{2}(Q^{\prime +}(0))\phi _{4}(\tfrac{S_{n-m-l}^{\prime
}}{a_{n}});L_{r-l}^{\prime }\geq 0,\tau _{n-m-l}^{\prime }=n-m-l]  \notag \\
& =\mathbf{E}[\phi _{2}(Q^{\prime +}(0));A_{0}]\mathbf{E}[\phi
_{4}(\tfrac{S_{N-R}^{\prime \prime }}{a_{n}});A_{1}]-\mathbf{E}[\phi _{2}(Q^{\prime
+}(0))\phi _{4}(\tfrac{S_{N-R}^{\prime \prime }}{a_{n}});A_{0}\cap A_{1}\cap
\overline{B}].  \notag \\
&  \label{iraz}
\end{align}%
Applying Lemma \ref{vmestoLem2.5} we see that, as $n\rightarrow \infty $
\begin{align}
& \mathbf{E}[\phi _{2}(Q^{\prime +}(0))\phi _{4}(\tfrac{S_{N-R}^{\prime
\prime }}{a_{n}});A_{0}\cap A_{1}\cap \overline{B}]  \notag \\
& \leq c\mathbf{P}\left( A_{1}\cap \overline{B}\right) =o\left(
\mathbf{P}\left( \tau _{N}=N\right) \right) =o\left( \mathbf{P}\left( \tau
_{n}=n\right) \right) ,  \label{idva}
\end{align}%
where $c>0.$ In much the same way as in proving of (\ref{perv}) it follows
from Lemma 2.3 in \cite{agkv} that \begin{equation}
\frac{\mathbf{E}[\phi _{2}(Q^{\prime +}(0));A_{0}]}{\mathbf{P}(L_{R}^{\prime
}\geq 0)}\rightarrow \mathbf{E}^{+}[\phi _{2}(Q^{+}(0))],\;n\rightarrow
\infty ,  \label{itri}
\end{equation}\begin{equation}
\frac{\mathbf{E}[\phi _{4}(\tfrac{S_{N-R}^{\prime \prime
}}{a_{n}});A_{1}]}{\mathbf{P}(\tau _{N-R}^{\prime \prime }=N-R)}\sim \frac{\mathbf{E}[\phi
_{4}(\tfrac{S_{n}}{a_{n}});\tau _{n}=n]}{\mathbf{P}(\tau _{n}=n)}\rightarrow
\mathbf{E}[\phi _{4}(L_{1}^{-})],\;n\rightarrow \infty .  \label{i4}
\end{equation}Relations (\ref{iraz})--(\ref{i4}) yield \begin{align}
& \frac{\mathbf{E}[\phi _{2}(Q^{\prime +}(0))\phi
_{4}(\tfrac{S_{n-m-l}^{\prime }}{a_{n}});L_{r-l}^{\prime }\geq 0,\tau _{n-m-l}^{\prime
}=n-m-l]}{\mathbf{P}(L_{R}^{\prime }\geq 0)\mathbf{P}(\tau _{n-m}=n-m)}
\notag \\
& \rightarrow \mathbf{E}^{+}[\phi _{2}(Q^{\prime +}(0))]\mathbf{E}[\phi
_{4}(L_{1}^{-})],\;n\rightarrow \infty .  \label{zavtra}
\end{align}%
It follows from Assumption A2 (see \cite{bi}), that $\tau _{n}/n$ converges
in distribution to a Beta-distribution with a density, which we denote by
$g(t)$. Since $\mathbf{P}(M_{l_{n}}<0)\mathbf{P}(L_{r-l_{n}}\geq
0)=\mathbf{P}(\tau _{r}=l_{n})$, it follows from \eqref{Dec1} -- (\ref{perv}) and\
(\ref{zavtra}) that \begin{align*}
\mathbf{E}\big[& \phi _{1}(Q^{-}(\tau _{r}))\phi _{2}(Q^{+}(\tau _{r}))\phi
_{3}(\tfrac{S_{\tau _{r}}}{a_{r}})\phi _{4}(\tfrac{S_{n-m}-S_{\tau
_{r}}}{a_{n}})|\tau _{n-m}=n-m\big] \\
& \quad \rightarrow \mathbf{E}^{-}\big[\phi
_{1}(Q_{1},...,Q_{k})\big]\mathbf{E}^{+}\big[\phi _{2}(Q_{1},...,Q_{k})\big] \\
& \qquad \qquad \mbox{}\times \mathbf{E}[\phi
_{4}(L_{1}^{-})]\,\int_{0}^{1}\mathbf{E}[\phi _{3}(t^{\frac{1}{\alpha }}L_{1}^{-})]g(t)\,dt\ .
\end{align*}

The lemma is proved.

\textbf{Remark. }Note that despite on the seeming similarity\ \ \ of the
formulations of Lemma \ref{leQ} and Lemma 2.6 in \cite{ABKVINTER} the
statements proved in it are essentially different, since the measure $\mu $
appearing as the limiting in Lemma \ref{leQ} is essentially different from
the corresponding limiting measure $\mu $, arising in Lemma 2.6 in
\cite{ABKVINTER}.

Set

\bigskip \begin{equation*}
\eta _{k}:=\sum_{y=0}^{\infty }y(y-1)Q_{k}(y)\Big/\Big(\sum_{y=0}^{\infty
}yQ_{k}(y)\Big)^{2}\ ,\quad k\geq 1\ .
\end{equation*}

$\ $

\begin{lemma}
\label{le2} (\cite{agkv}, \cite{ABKVINTER}) Assume Assumptions A1 -- A3.
Then for all $x\geq 0$ \begin{equation*}
\sum_{k=0}^{\infty }\eta _{k+1}e^{-S_{k}}\ <\ \infty \qquad
\mathbf{P}_{x}^{+}\text{ -a.s.}
\end{equation*}and for all $x\leq 0$ \begin{equation*}
\sum_{k=1}^{\infty }\eta _{k}e^{S_{k}}\ <\ \infty \qquad
\mathbf{P}_{x}^{-}\text{ -a.s.}
\end{equation*}
\end{lemma}

\section{Trees with steam}

We need the following construction considered in \cite{ABKVINTER} which is
in spirit to the models from \cite{chau91}, \cite{ka77}. For $n=0,1,\ldots
,\infty $ let $\mathcal{T}_{n}$ be the set of all ordered rooted trees of
height exactly $n$, whose edges \ are directe from the root. The precise
definition of ordered rooted trees was given by Neven \cite{ne}. Let
$\mathcal{T}_{\geq n}=\mathcal{T}_{n}\cup \mathcal{T}_{n+1}\cup \cdots \cup
\mathcal{T}_{\infty }$ be the set of ordered rooted trees of at least height
$n$, whose edges are directed from the root. Denote by $[t]_{n}\in
\mathcal{T}_{n}$\ \ the tree obtained from the tree \ $t\in \mathcal{T}_{\geq n}$ $\
$by eliminating from the tree $t$ all edges and nodes of a height exceeding\
$n$. We will say that the tree $[t]_{n}$ is obtained by \ pruning a tree $t$
to the tree $[t]_{n}$ $\ $on the level $n$. \

For $n=0,1,\ldots ,\infty $ a tree with a stem of height $n$, or shortly an
\emph{o-tree} of height $n$, is a pair \begin{equation*}
\mathsf{t}=(t,k_{0}k_{1}\ldots k_{n})\ ,
\end{equation*}%
where $t\in \mathcal{T}_{\geq n}$, $k_{0},\ldots ,k_{n}$ \ are nodes in $t$
such that $k_{0}$ is the root, and nodes $k_{i-1}$ and $k_{i},$ $i=1,...,n,$
are connected by the edge directed from $k_{i-1}$ to $k_{i}$, i.e. the node
$k_{i}$ belongs to the generation $i$. We call $k_{0}\ldots k_{n}$ \ the stem
within o-tree $\mathsf{t}$ ( which, evidently, is determined by $k_{n}).$
Let $\mathcal{T}_{n}^{\prime }$ be the set of all o-trees with the trests of
height $n$.

Pruning an o-tree $\mathsf{t}=(t,k_{0}k_{1}\ldots k_{n})$ of the heigth $n$
at level $m\leq n,$ we \ obtain the o-tree \begin{equation*}
\lbrack \mathsf{t}]_{m}=([t]_{m},k_{0}\ldots k_{m})\
\end{equation*}of the heigth $m.$

Every tree $t\in \mathcal{T}_{\geq n}$ includes a unique o-tree
\begin{equation*}
\langle t\rangle _{n}=([t]_{n},k_{0}(t)\ldots k_{n}(t))
\end{equation*}%
of height $n$, where $k_{0}(t)\ldots k_{n}(t)$ is the\emph{\ leftmost }stem,
which can be fitted into $[t]_{n}$. \newline

In the sequel it will be convinuent to assume that there is a particle in
every node of the tree $t.$The particle located at the node $a$ is the child
of the particle in the node $b\in t$, \ if the nodes $a$\ and $b$ are
connected by the edge directed from $a$ to $b$.

Let $\pi =(q_{1},q_{2},\ldots )$ be a fixed environment. Define the discrete
distribution $\tilde{q}_{i}$\ \ by its weights \begin{equation*}
\tilde{q}_{i}(y)=\tfrac{1}{m(q_{i})}yq_{i}(y)\ ,\quad y=0,1,\ldots .
\end{equation*}%
Then an \emph{oLLP-tree\ }(Lyons-Pemantle-Peres o-tree) corresponding to the
disrtibution $(q_{1},q_{2},\ldots )$ is the random o-tree
$\tilde{\mathsf{T}}=(\tilde{T},\tilde{K}_{0}\tilde{K}_{1}\ldots )$ with values in
$\mathcal{T}_{\infty }^{\prime },$ satisfying the following properties:

Given $\Pi =(q_{1},q_{2},\ldots )$

\begin{itemize}
\item The offspring numbers of all particles are independent random
variables.

\item The offspring number of $\tilde{K}_{i-1}$ has disribution
$\tilde{q}_{i}$ and the offspring number of any other particle in generation $i-1$ has
disribution $q_{i}$.

\item The node $\tilde{K}_{i}$ is uniformly distributed among all children
of $\tilde{K}_{i-1}$.
\end{itemize}

In other words the particles being in the noodes of infinite stem reproduce
according to the distribution $(\tilde{q}_{1},\tilde{q}_{2},\ldots )$, and
all other particles reproduce according to the distribution
$(q_{1},q_{2},\ldots )$.

Recall some useful properties of oLLP-trees.\ Let $\tilde{Z}_{n}$ be the
population size of the oLLP-tree in generation $n$.

We have the following statement.

\begin{lemma}
\label{le41} (\cite{ABKVINTER}) If Assumptions A1 to A3 are valid, then as
$n\rightarrow \infty $ \begin{equation*}
e^{-S_{n}}\tilde{Z}_{n}\rightarrow W^{+}\qquad \mathbf{P}^{+}\text{-a.s.,}
\end{equation*}%
where the random variable $W^{+}$ is such that $W^{+}>0$
$\mathbf{P}^{+}$-a.s.
\end{lemma}

We use the representation \begin{equation*}
\tilde{Z}_{n}=1+\sum_{i=0}^{n-1}\tilde{Z}_{n}^{i},
\end{equation*}%
where $\tilde{Z}_{n}^{i}$ is the number of particles in generation $n$ other
than $\tilde{K}_{n}$, which are descent of the particle from
$\tilde{K}_{i}$, but not of the particle from $\tilde{K}_{i+1}$.--
число вершин в поколении $n,$ отличных от вершины $\tilde{K}_{n}$,
которые являются потомками частицы вершины $\tilde{K}_{i}$, но
не частицы вершины $\tilde{K}_{i+1}$. Note that
$\mathbf{E}[\tilde{Z}_{i+1}^{i}\mid \Pi ]=\sum_{y}y\tilde{Q}_{i+1}(y)-1=e^{X_{i+1}}\eta
_{i+1}$ and a.s. \begin{equation}
\mathbf{E}[\tilde{Z}_{n}^{i}\mid \Pi
]=e^{S_{n}-S_{i+1}}\mathbf{E}[\tilde{Z}_{i+1}^{i}\mid \Pi ]=\eta _{i+1}e^{S_{n}-S_{i}}\ .  \label{4.9is}
\end{equation}

Clearly, that given the environment, the sequence
$e^{-S_{n}}\sum_{i=k}^{n-1}\tilde{Z}_{n}^{i}$ is for $n>k$ a non-negative submartingale. From here and
Doob's inequality it follows that for every $\varepsilon \in (0,1)$
\begin{equation}
\mathbf{P}\Big(\max_{k<m\leq
n}e^{-S_{m}}\sum_{i=k}^{m-1}\tilde{Z}_{m}^{i}\geq \varepsilon \ \Big|\ \Pi \Big)\leq \frac{1}{\varepsilon
}\sum_{i=k}^{n-1}e^{-S_{n}}\mathbf{E}[\tilde{Z}_{n}^{i}\mid \Pi ]\leq
\frac{1}{\varepsilon }\sum_{i\geq k}\eta _{i+1}e^{-S_{i}}  \label{4.9aa}
\end{equation}and \begin{equation}
\mathbf{P}^{+}\Big(\sup_{m>k}e^{-S_{m}}\sum_{i=k}^{m-1}\tilde{Z}_{m}^{i}\geq
\varepsilon \Big)\leq \frac{1}{\varepsilon }\mathbf{E}^{+}\Big[1\wedge
\sum_{i\geq k}\mathbf{\eta }_{i+1}e^{-S_{i}}\Big]\ .  \label{4.9bb}
\end{equation}It follows from \ref{le2} that for sufficiently large $k$
\begin{equation*}
\mathbf{P}^{+}\Big(\sup_{m>k}e^{-S_{m}}\sum_{i=k}^{m-1}\tilde{Z}_{m}^{i}\geq
\varepsilon \Big)\leq \varepsilon \ .
\end{equation*}

To approximate the conditional distribution of the considered branching
process $Z=\left( Z_{0},Z_{1},...\right) $ given $Z_{n}>0$ we use its
approximation by oLLP-trees. Let $T$ \ be the tree corresponding to the
considered branching process $Z$ in random environment $\Pi .$

We need the following statement.

\begin{theorem}
\label{trest} (\cite{ABKVINTER}) Let $0\leq d_{n}<n$ be a sequence of
natural numbers with $d_{n}\rightarrow \infty $ as $n\rightarrow \infty $,
and $Y_{n}$ be uniformly bounded random variables of the form $Y_{n}=\varphi
(Q_{1},\ldots ,Q_{n-d_{n}}),$ and let $B_{n}\subset
\mathcal{T}_{n-d_{n}}^{\prime }$, $n\geq 1$. If \ Assumptions A1 -- A3 are valid and
for some $\ell \geq 0$ \begin{equation*}
\mathbf{E}\big[Y_{n};[\tilde{\mathsf{T}}]_{n-d_{n}}\in B_{n}\ \big|\ \tau
_{n-m}=n-m\big]\rightarrow \ell
\end{equation*}for all $m\geq 0$, then \begin{equation*}
\mathbb{E}\big[Y_{n};[\langle T\rangle _{n}]_{n-d_{n}}\in B_{n}\ \big|\
Z_{n}>0\big]\rightarrow \ell \ .
\end{equation*}%
Here the set $B_{n}$ may be random, depending only on the environment $\Pi $.
\end{theorem}

\noindent

\noindent

\section{Proof of Theorem 1.4}

Let $\tilde{\mathsf{T}}$ be an oLLP---tree. Recall that for $i<j$
$\tilde{Z}_{j}^{i}$ is  the number of the particles in generation $j$ other than the
particle in $\tilde{K}_{j}$, which descent from a particle located in
$\tilde{K}_{i}$ but not from a particle located in $\tilde{K}_{i+1}$. For
convenience we put $\tilde{Z}_{j}^{i}=0$ for $i\geq j$.

Recall also that $r=r_{n}\rightarrow \infty ,n\rightarrow \infty
,r_{n}=o\left( n\right) ,\tau _{r}=\tau _{r_{n}}.$

\begin{lemma}
\label{le51} For every $\varepsilon >0$ there is a natural number $a$ such
that for any natural numbers $m$ and $\varsigma \in \lbrack \tau _{r},r]$
for sufficiently large $n$ (depending on $\varepsilon ,a$ and $m$)
\begin{equation*}
\mathbf{P}\Big(\sum_{i:|i-\tau _{r}|\geq a}\frac{\tilde{Z}_{\varsigma
}^{i}}{e^{S_{\varsigma }-S_{\tau _{r}}}}\geq \varepsilon \ \Big|\ \tau
_{n-m}=n-m\Big)\leq \varepsilon \ ,
\end{equation*}%
here $\varsigma $ may be random, depending only on the random environment
$\Pi $.
\end{lemma}

\textbf{Proof.} For $0\leq j<n$ set \begin{equation*}
L_{j,n}=\min (S_{k+1},...,S_{n})-S_{n},\;L_{n,n}=0.
\end{equation*}%
from Markov inequality and (\ref{4.9is}) we have for $0<\varepsilon \leq 1$
and $m\leq n-r_{n}$ \begin{align*}
\varepsilon \mathbf{P}\Big(\sum_{|i-\tau _{r}|\geq
a}\frac{\tilde{Z}_{\varsigma }^{i}}{e^{S_{\varsigma }-S_{\tau _{r}}}}& \geq \varepsilon ;\tau
_{n-m}={n-m}\Big) \\
& \leq \mathbf{E}\Big[1\wedge \sum_{i\leq \varsigma ,|i-\tau _{r}|\geq
a}\eta _{i+1}e^{S_{\tau _{r}}-S_{i}};\tau _{n-m}=n-m\Big]\  \\
& \leq \sum_{j\leq r}\mathbf{E}\Big[1\wedge \sum_{i\leq \varsigma ,|i-j|\geq
a}\eta _{i+1}e^{S_{j}-S_{i}};\tau _{j}=j,L_{j,r}\geq 0\Big] \\
& \qquad \qquad \qquad \times \mathbf{P}\big(\tau _{n-r-m}=n-r-m\big)\ .
\end{align*}Further, \begin{align*}
\sum_{j\leq r}& \mathbf{E}\Big[1\wedge \sum_{i\leq \varsigma ,|i-j|\geq
a}\eta _{i+1}e^{S_{j}-S_{i}};\tau _{j}=j,L_{j,r}\geq 0\Big]= \\
& =\sum_{j\leq r}\mathbf{E}\Big[1\wedge \sum_{i=0}^{j-a}\eta
_{i+1}e^{S_{j}-S_{i}};\tau _{j}=j\Big]\mathbf{P}(L_{r-j}\geq 0) \\
& \qquad \qquad \mbox{}+\sum_{j\leq r}\mathbf{P}(\tau
_{j}=j)\mathbf{E}\Big[1\wedge \sum_{i=j+a}^{\varsigma }\eta _{i+1}e^{S_{j}-S_{i}};L_{j,r}\geq
0\Big]\ .
\end{align*}Using duality, we obtain \begin{align*}
\sum_{j\leq r}& \mathbf{E}\Big[1\wedge \sum_{i\leq \varsigma ,|i-j|\geq
a}\eta _{i+1}e^{S_{j}-S_{i}};\tau _{j}=j,L_{j,r}\geq 0\Big]= \\
& \leq \sum_{a\leq j\leq r}\mathbf{E}\Big[1\wedge \sum_{i=a}^{j-a}\eta
_{i}e^{S_{i}};M_{j}<0\Big]\mathbf{P}(L_{r-j}\geq 0) \\
& \qquad \qquad \mbox{}+\sum_{a\leq k\leq r}\mathbf{P}(\tau
_{r}-k=r-k)\mathbf{E}\Big[1\wedge \sum_{i=a}^{k}\eta _{i+1}e^{-S_{i}};L_{k}\geq 0\Big]\
.
\end{align*}%
By Lemmas \ref{limitEminus} and \ref{le2} we may choose $a$ so large that
for a given $\delta >0$ and all $j,k>a$ \begin{align*}
\mathbf{E}\Big[1\wedge \sum_{i=a}^{j}\eta _{i}e^{S_{i}};M_{j}<0\Big]& \leq
\delta \mathbf{P}(M_{j}<0), \\
\mathbf{E}\Big[1\wedge \sum_{i=a}^{k}\eta _{i+1}e^{-S_{i}};L_{k}\geq 0\Big]&
\leq \delta \mathbf{P}(L_{k}\geq 0).
\end{align*}Duality yields \begin{align*}
\sum_{j\leq r}& \mathbf{E}\Big[1\wedge \sum_{i\leq \varsigma ,|i-j|\geq
a}\eta _{i+1}e^{S_{j}-S_{i}};\tau _{j}=j,L_{j,r}\geq 0\Big] \\
& \leq \delta \sum_{a\leq j\leq r}\mathbf{P}(\tau
_{j}=j)\mathbf{P}(L_{r-j}\geq 0) \\
& \qquad \qquad \mbox{}+\delta \sum_{a\leq k\leq r}\mathbf{P}(\tau
_{r-k}=r-k)\mathbf{P}(L_{k}\geq 0)\leq 2\delta \end{align*}and \begin{align*}
\mathbf{P}\Big(\sum_{|i-\tau _{r}|\geq a}\frac{\tilde{Z}_{\varsigma
}^{i}}{e^{S_{\varsigma }-S_{\tau _{r}}}}& \geq \varepsilon ;\tau _{n-m}=n-m\Big) \\
& \leq \frac{2\delta }{\varepsilon }\mathbf{P}\big(\tau _{n-r-m}=n-r-m\big)\
.
\end{align*}%
Since $\mathbf{P}(\tau _{n}=n)$ is regularly varying and $r=o\left( n\right)
,$ $n\rightarrow \infty $, the right-hand side is bounded by $\varepsilon
\mathbf{P}(\tau _{n}=n)$, if $\delta $ is chosen small enough. The lemma is
proved.

We now come back to the proof of the first part of Theorem \ref{theomain}.
Let, as before, $\tilde{Z}_{j}$ be the number of nodes in generation $j$ of
the oLPP-tree $\tilde{\mathsf{T}}$, i.e., \begin{equation*}
\tilde{Z}_{j}=1+\sum_{k=0}^{j-1}\tilde{Z}_{j}^{k}\ .
\end{equation*}%
Fix an $\varepsilon >0.$ In view of the preceding lemma with $\varsigma
=\tau _{r}$ there is a natural number $a$ such that given $\tau _{n-m}=n-m$
the probability is at least $1-\varepsilon $ that the event \begin{equation*}
\tilde{Z}_{\tau _{r}}=1+\sum_{|k-\tau _{r}|\leq a}\tilde{Z}_{\tau
_{r}}^{k}=1+\sum_{k=\tau _{r}-a}^{\tau _{r}}\tilde{Z}_{\tau _{r}}^{k}
\end{equation*}%
holds. Note that given the environment $\Pi $ the distribution of
\begin{equation*}
1+\sum_{k=\tau _{r}-a}^{\tau _{r}}\tilde{Z}_{\tau _{r}}^{k}
\end{equation*}%
only depends on $(Q_{\tau _{r}-a},\ldots ,Q_{\tau _{r}})$. \ By Lemma
\ref{leQ} given $\tau _{n-m}=n-m$ the vector $(Q_{\tau _{r}-a},\ldots ,Q_{\tau
_{r}})$ \ converges in distribution to a limit.

These observations hold for every $\varepsilon >0$. Therefore we may
summarize our arguments \ as follows: For all $m\geq 1$ \begin{equation*}
\big(\tilde{Z}_{\tau _{r}}\mid \tau _{n-m}=n-m\big)\overset{d}{\rightarrow
}\varkappa _{1}\ ,
\end{equation*}%
where the random variable $\varkappa _{1}$ has the properties as clamed in
Theorem \ref{theomain}. \ Now Theorem \ref{trest} gives the claim of the
first part of Theorem \ref{theomain}.

Now we pass to the proof of the second part of Theorem \ref{theomain}.\ Let
for fixed $a$ \begin{equation*}
\hat{Z}_{a,k}=\sum_{i:|i-\tau _{r}|\leq a}\tilde{Z}_{k}^{i}
\end{equation*}and \begin{equation*}
\alpha _{a,r}=e^{S_{\tau _{r}}-S_{r}}\hat{Z}_{a,r}\ ,\ \beta
_{a,r}=e^{S_{\tau _{r}}-S_{\tau _{r}+a}}\hat{Z}_{a,\tau _{r}+a}\ .
\end{equation*}We need the following statement.

\begin{lemma}
\label{le52} Let $m\geq 1$ and $\varepsilon \in \left( 0,1\right) $. Then,
if $a$ is sufficiently large \begin{equation}
\limsup_{n\rightarrow \infty }\mathbf{P}(|\alpha _{a,r}-\beta
_{a,r}|>\varepsilon \mid \tau _{n-m}=n-m)\leq \varepsilon \ .  \label{nuii}
\end{equation}
\end{lemma}

\textbf{Proof.} In virtue of Markov inequality and (\ref{4.9is}),
(\ref{4.9aa}) and (\ref{4.9bb}) \begin{align*}
\mathbf{P}(\beta _{a,r}& >d\mid \tau _{n-m}=n-m) \\
& \leq \mathbf{P}(e^{S_{\tau _{r}}-S_{\tau
_{r}+a}}\mathbf{E}[\hat{Z}_{a,\tau _{r}+a}\mid \Pi ]>\sqrt{d}\mid \tau _{n-m}=n-m)+\frac{1}{\sqrt{d}}
\\
& \leq \mathbf{P}\Big(\sum_{i:|i-\tau _{r}|\leq a}\eta _{i+1}e^{S_{\tau
_{r}}-S_{i}}>\sqrt{d}\ \Big|\ \tau _{n-m}=n-m\Big)+\frac{1}{\sqrt{d}}\ .
\end{align*}%
It follows \ from Lemma \ref{leQ} that the sum under the sign of probability
converges in distribution as $n\rightarrow \infty $ and \begin{align*}
\limsup_{n\rightarrow \infty }\mathbf{P}(\beta _{a,n}& >d\mid \tau
_{n-m}=n-m) \\
& \leq \mathbf{P}^{-}\Big(\sum_{i\geq 1}\eta _{i}e^{S_{i}}\geq
\frac{\sqrt{d}}{2}\Big)+\mathbf{P}^{+}\Big(\sum_{i\geq 0}\eta _{i+1}e^{-S_{i}}\geq
\frac{\sqrt{d}}{2}\Big)+\frac{1}{\sqrt{d}}\ .
\end{align*}%
By Lemma \ref{le2} there exists $d<\infty $ such that for all $a>0$
\begin{equation*}
\limsup_{n\rightarrow \infty }\mathbf{P}(\beta _{a,r}>d\mid \tau
_{n-m}=n-m)<\varepsilon /2\ .
\end{equation*}

Further, since $r=o\left( n\right) $ and $\mathbf{P}(\tau _{n}=n)$ is
regulary varying at infinity we have for sufficiently large $n$
\begin{eqnarray}
\mathbf{P}(\tau _{r} &\in &[r-a,r]\mid \tau
_{n-m}=n-m)=\frac{\mathbf{P}(\tau _{r}\in \lbrack r-a,r],\tau _{n-m}=n-m)}{\mathbf{P}(\tau _{n-m}=n-m)}
\notag \\
&\leq &\sum_{k=r-a}^{r}\frac{\mathbf{P}(\tau _{r}=k)\mathbf{P}(\tau
_{n-m-r}=n-m-r)}{\mathbf{P}(\tau _{n-m}=n-m)}  \notag \\
&\leq &2\sum_{k=r-a}^{r}\mathbf{P}(\tau _{r}=k)=2\mathbf{P}(\tau _{r}\in
\lbrack r-a,r]).  \label{zebra}
\end{eqnarray}

Since the distribution of the ratio $\tau _{n}/n$ converges as $n\rightarrow
\infty $ to a Beta distribution (\ref{zebra}) implies \begin{equation*}
\mathbf{P}(\tau _{r}+a\geq r\mid \tau _{n-m}=n-m)\rightarrow
0,\;n\rightarrow \infty .
\end{equation*}

Therefore  \begin{align}
\mathbf{P}(& |\beta _{a,r}-\alpha _{a,r}|>\varepsilon \mid \tau _{n-m}=n-m)
\notag \\
& \leq \frac{\varepsilon }{2}+\mathbf{P}\big(|\alpha _{a,r}-\beta
_{a,r}|>\varepsilon ,\beta _{a,r}\leq d,\tau _{r}+a\leq r\mid \tau
_{n-m}=n-m\big)\   \label{help0301}
\end{align}%
for sufficiently large $n.$ Note, that given $\Pi $, $\hat{Z}_{a,\tau _{r}+a}
$ and $\tau _{r}+a\leq r$, the process $\hat{Z}_{a,k}$, $k\geq \tau _{r}+a$
is a branching process in varying environment. Therefore $\mathbf{E}[\alpha
_{a,r}\mid \Pi ,\hat{Z}_{a,\tau _{r}+a}]=\beta _{a,r}$ a.s. and
\begin{equation*}
\frac{\mathbf{Var}(Z_{n}\mid Z_{0}=z,\Pi )}{\mathbf{E}[Z_{n}\mid Z_{0}=1,\Pi
]^{2}}=z\Big(e^{-S_{n}}+\sum_{i=0}^{n-1}\eta _{i+1}e^{-S_{i}}-1\Big)\ .
\end{equation*}

Inserting this estimate into (\ref{help0301}), we have

\bigskip \begin{align}
\mathbf{P}(& |\beta _{a,r}-\alpha _{a,r}|>\varepsilon ;\tau _{n-m}=n-m)
\notag \\
& \leq \frac{\varepsilon }{2}\mathbf{P}(\tau
_{n-m}=n-m)+\frac{d}{\varepsilon ^{2}}\mathbf{E}\Big[1\wedge \Big(e^{-(S_{r}-S_{\tau _{r}})}
\notag \\
& \qquad \qquad \mbox{}+\sum_{i=\tau _{r}+a}^{r}\eta
_{i+1}e^{-(S_{i}-S_{\tau _{r}})}\Big);\tau _{r}+a\leq r,\tau
_{n-m}=n-m\Big].\qquad  \label{vott2}
\end{align}

\bigskip Thus to finish the proof of Lemma \ it is sufficient to show that
the term \begin{eqnarray*}
&&\mathbf{E}\left[ 1\wedge \Big(e^{-(S_{r}-S_{\tau _{r}})}+\sum_{i=\tau
_{nt}+a}^{r}\eta _{i+1}e^{-(S_{i}-S_{\tau _{r}})}\Big);\tau _{r}+a\leq
r,\tau _{n-m}=n-m\right] \\
&=&\sum_{j\leq r-a}\mathbf{E}\left[ 1\wedge
\Big(e^{-(S_{r}-S_{j})}+\sum_{i=j+a}^{r}\eta _{i+1}e^{-(S_{i}-S_{j})}\Big);\tau _{r}=j,\tau
_{n-m}=n-m\right]
\end{eqnarray*}may be made arbitrary small.

Since for $j\leq r=o(n)$\begin{equation*}
\sup_{0\leq j\leq r}\frac{\mathbf{P}\left( \tau _{n-m-j}=n-m-j\right)
}{\mathbf{P}\left( \tau _{n-m}=n-m\right) }\rightarrow 1
\end{equation*}%
as $n\rightarrow \infty ,$ we have for all sufficiently large $n$ and $r=o(n)
$\begin{eqnarray}
&&\mathbf{E}\left[ 1\wedge \Big(e^{-(S_{r}-S_{j})}+\sum_{i=j+a}^{r}\eta
_{i+1}e^{-(S_{i}-S_{j})}\Big);\tau _{r}=j,\tau _{n-m}=n-m\right]   \notag \\
&\leq &\mathbf{E}\left[ 1\wedge \Big(e^{-(S_{r}-S_{j})}+\sum_{i=j+a}^{r}\eta
_{i+1}e^{-(S_{i}-S_{j})}\Big);\tau _{r}=j\right] \mathbf{P}\left( \tau
_{n-m-j}=n-m-j\right)   \notag \\
&\leq &2\mathbf{E}\left[ 1\wedge
\Big(e^{-(S_{r}-S_{j})}+\sum_{i=j+a}^{r}\eta _{i+1}e^{-(S_{i}-S_{j})}\Big);\tau _{r}=j\right] \mathbf{P}\left( \tau
_{n-m}=n-m\right) .  \notag \\
&&  \label{dop1}
\end{eqnarray}Further,\begin{eqnarray}
&&\mathbf{E}\left[ 1\wedge \Big(e^{-(S_{r}-S_{j})}+\sum_{i=j+a}^{r}\eta
_{i+1}e^{-(S_{i}-S_{j})}\Big);\tau _{r}=j\right]   \notag \\
&=&\mathbf{P}\left( \tau _{j}=j\right) \mathbf{E}\left[ 1\wedge
\Big(e^{-(S_{r}-S_{j})}+\sum_{i=j+a}^{r}\eta
_{i+1}e^{-(S_{i}-S_{j})}\Big);L_{j,r}\geq 0\right]   \notag \\
&\leq &\mathbf{P}\left( \tau _{j}=j\right) \mathbf{P}(L_{r-j}\geq 0)  \notag
\\
&&\times \left( \mathbf{E}\left[ e^{-S_{r-j}}|L_{r-j}\geq 0\right]
+\mathbf{E}\left[ 1\wedge \left( \sum_{i=j+a}^{r}\eta _{i-j+1}e^{-S_{i-j}}\right)
|L_{r-j}\geq 0\right] \right) .  \notag \\
&&  \label{dobav1}
\end{eqnarray}%
Note, that by Lemma \ref{limitEminus} with $U_{n}=e^{-S_{n}}$ we have for
sufficiently large $a$ \ \begin{equation}
\sup_{0\leq j\leq r-a}\mathbf{E}\left[ e^{-S_{r-j}}|L_{r-j}\geq 0\right]
\leq \sup_{a\leq k\leq \infty }\mathbf{E}\left[ e^{-S_{k}}|L_{k}\geq 0\right]
\leq \varepsilon ^{3}/\left( 8d\right) .  \label{dobav}
\end{equation}Finally, \begin{equation*}
\mathbf{E}\left[ 1\wedge \left( \sum_{i=j+a}^{r}\eta
_{i-j+1}e^{-S_{i-j}}\right) |L_{r-j}\geq 0\right] =\mathbf{E}\left[ 1\wedge
\left( \sum_{k=a}^{r-j}\eta _{k+1}e^{-S_{k}}\right) |L_{r-j}\geq 0\right]
\end{equation*}and \begin{equation*}
\lim_{r-j\rightarrow \infty }1\wedge \left( \sum_{k=a}^{r-j}\eta
_{k+1}e^{-S_{k}}\right) =1\wedge \left( \sum_{k=a}^{\infty }\eta
_{k+1}e^{-S_{k}}\right)
\end{equation*}$\mathbf{P}^{+}-$a.s. Hence, applying Lemma \ref{limitEminus} we obtain

\begin{equation*}
\lim_{r-j\rightarrow \infty }\mathbf{E}\left[ 1\wedge \left(
\sum_{k=a}^{r-j}\eta _{k+1}e^{-S_{k}}\right) |L_{r-j}\geq 0\right]
=\mathbf{E}^{+}\left[ 1\wedge \left( \sum_{k=a}^{\infty }\eta _{k+1}e^{-S_{k}}\right)
\right] .
\end{equation*}Since\begin{equation*}
\sum_{k=a}^{\infty }\eta _{k+1}e^{-S_{k}}<\infty
\;\mathbf{P}^{+}-\text{a.s.},
\end{equation*}it follows that\begin{equation*}
\lim_{a\rightarrow \infty }1\wedge \left( \sum_{k=a}^{\infty }\eta
_{k+1}e^{-S_{k}}\right) =0\;\mathbf{P}^{+}-\text{a.s.}
\end{equation*}Hence, by the dominated convergence theorem we get
\begin{equation*}
\lim_{a\rightarrow \infty }\mathbf{E}^{+}\left[ 1\wedge \left(
\sum_{k=a}^{\infty }\eta _{k+1}e^{-S_{k}}\right) \right]
=\mathbf{E}^{+}\left[ \lim_{a\rightarrow \infty }\left( 1\wedge \left( \sum_{k=a}^{\infty
}\eta _{k+1}e^{-S_{k}}\right) \right) \right] =0.
\end{equation*}Thus,\begin{equation*}
\lim_{a\rightarrow \infty }\lim_{r-j\rightarrow \infty }\mathbf{E}\left[
1\wedge \left( \sum_{k=a}^{r-j}\eta _{k+1}e^{-S_{k}}\right) |L_{r-j}\geq
0\right] =0.
\end{equation*}

Therefore, given \ $r-j\geq 2a$ we have for sufficiently large $a$
\begin{eqnarray}
\mathbf{E}\left[ 1\wedge \left( \sum_{k=a}^{r-j}\eta _{k+1}e^{-S_{k}}\right)
|L_{r-j}\geq 0\right] &\leq &2\mathbf{E}^{+}\left[ 1\wedge \left(
\sum_{k=a}^{\infty }\eta _{k+1}e^{-S_{k}}\right) \right] \leq \varepsilon
^{3}/\left( 8d\right) .  \notag \\
&&  \label{dop2}
\end{eqnarray}Thus, it remains to consider the case \begin{equation*}
r-2a\leq j\leq r-a.
\end{equation*}%
Since the random variable $\tau _{n}/n$ converges in distribution as
$n\rightarrow \infty $ to Beta distribution we have for sufficiently large $r$
\begin{eqnarray*}
&&\sum_{r-2a\leq j\leq r-a}\mathbf{E}\left[ 1\wedge
\Big(e^{-(S_{r}-S_{j})}+\sum_{i=j+a}^{r}\eta _{i+1}e^{-(S_{i}-S_{j})}\Big);\tau _{r}=j\right] \\
&\leq &\sum_{r-2a\leq j\leq r-a}\mathbf{P}\left( \tau _{r}=j\right)
=\mathbf{P}\left( r-2a\leq \tau _{r}\leq r-a\right) \leq \varepsilon ^{3}/\left(
4d\right) .
\end{eqnarray*}This gives the needed result \begin{eqnarray}
&&\sum_{j\leq r-a}\mathbf{E}\left[ 1\wedge
\Big(e^{-(S_{r}-S_{j})}+\sum_{i=j+a}^{r}\eta _{i+1}e^{-(S_{i}-S_{j})}\Big);\tau _{r}=j\right]  \notag
\\
&\leq &\varepsilon ^{3}\mathbf{P}\left( \tau _{r}\leq r-2a\right) /\left(
4d\right) +\varepsilon ^{3}/\left( 4d\right) \mathbf{\leq }\varepsilon
^{3}/2.  \label{vot1}
\end{eqnarray}

Relations (\ref{vott2}) -- (\ref{vot1}) yield (\ref{nuii}).

We are now ready to finish the proof of the second part of Theorem
\ref{theomain}. From $\tilde{Z}_{r}=1+\hat{Z}_{a,r}+\sum_{i:|i-\tau
_{r}|>a}\tilde{Z}_{r}^{i}$ we have \begin{align*}
\mathbf{P}\big(|e^{S_{\tau _{r}}-S_{r}}& \tilde{Z}_{r}-\beta _{a,r}|\geq
3\varepsilon \mid \tau _{n-m}=n-m) \\
& \leq \mathbf{P}(e^{S_{\tau _{r}}-S_{r}}\geq \varepsilon \mid \tau
_{n-m}=n-m) \\
& \mbox{}\quad +\mathbf{P}(|\alpha _{a,r}-\beta _{a,r}|\geq \varepsilon \mid
\tau _{n-m}=n-m) \\
& \mbox{}\quad +\mathbf{P}\Big(e^{S_{\tau _{r}}-S_{r}}\sum_{i:|i-\tau
_{r}|>a}\tilde{Z}_{r}^{i}\geq \varepsilon \ \Big|\ \tau _{n-m}=n-m\Big)\ .
\end{align*}

\bigskip Fix an $\varepsilon >0.$ Then\begin{eqnarray}
\mathbf{P}(e^{S_{\tau _{r}}-S_{r}} &\geq &\varepsilon \mid \tau _{n-m}=n-m)
\notag \\
&=&\frac{1}{\mathbf{P}(\tau _{n-m}=n-m)}\mathbf{P}\big(\tfrac{S_{\tau
_{r}}-S_{r}}{a_{r}}\geq \tfrac{\log \varepsilon }{a_{r}};\tau
_{n-m}=n-m\big).  \notag \\
&&  \label{a1}
\end{eqnarray}Clearly,\begin{eqnarray}
\mathbf{P}\big(\tfrac{S_{\tau _{r}}-S_{r}}{a_{r}} &\geq &\tfrac{\log
\varepsilon }{a_{r}};\tau _{n-m}=n-m\big)  \notag \\
&=&\sum\limits_{k=0}^{r}\mathbf{P}\big(\tfrac{S_{k}-S_{r}}{a_{r}}\geq
\tfrac{\log \varepsilon }{a_{r}};\tau _{r}=k;\tau _{n-m}=n-m\big).  \notag \\
&&  \label{a2}
\end{eqnarray}Let $S_{i}^{\prime }=S_{r+i}-S_{r},n\geq 0,$ and
\begin{equation*}
\tau _{n}^{\prime }=\min \left\{ l\in \left[ 0,n\right] :S_{l}^{\prime
}=\min_{0\leq i\leq n}S_{i}^{\prime }\right\} .
\end{equation*}For a fixed $k$ we write the following chain of relations
\begin{eqnarray}
\mathbf{P}\big(\tfrac{S_{k}-S_{r}}{a_{r}} &\geq &\tfrac{\log \varepsilon
}{a_{r}};\tau _{r}=k,\tau _{n-m}=n-m\big)  \notag \\
&\leq &\mathbf{P}\big(\tfrac{S_{k}-S_{r}}{a_{r}}\geq \tfrac{\log \varepsilon
}{a_{r}};\tau _{r}=k,\text{ }\tau _{n-m-r}^{\prime }=n-m-r\big)  \notag \\
&=&\mathbf{P}\big(\tfrac{S_{k}-S_{r}}{a_{r}}\geq \tfrac{\log \varepsilon
}{a_{r}};\tau _{r}=k\big)\mathbf{P}\left( \tau _{n-m-r}^{\prime }=n-m-r\right)
\notag \\
&=&\mathbf{P}\big(\tfrac{S_{k}-S_{r}}{a_{r}}\geq \tfrac{\log \varepsilon
}{a_{r}};\tau _{r}=k\big)\mathbf{P}\left( \tau _{n-m-r}=n-m-r\right) .
\label{a3}
\end{eqnarray}%
Since the probability $\mathbf{P}\left( \tau _{n}=n\right) =\mathbf{P}\left(
M_{n}<0\right) $ is regulary varying as $n\rightarrow \infty $\ and
$r=o\left( n\right) ,$ we have \begin{equation*}
\lim_{n\rightarrow \infty }\frac{\mathbf{P}\left( \tau _{n-m-r}=n-m-r\right)
}{\mathbf{P}\left( \tau _{n}=n\right) }=1.
\end{equation*}%
It follows now from (\ref{a1})--(\ref{a3}) that for large $n$
\begin{eqnarray*}
\mathbf{P}(e^{S_{\tau _{r}}-S_{r}} &\geq &\varepsilon \mid \tau _{n-m}=n-m)
\\
&\leq &2\sum\limits_{k=0}^{r}\mathbf{P}\big(\tfrac{S_{k}-S_{r}}{a_{r}}\geq
\tfrac{\log \varepsilon }{a_{r}};\tau _{r}=k\big) \\
&=&2\mathbf{P}\big(\tfrac{S_{\tau _{r}}-S_{r}}{a_{r}}\geq \tfrac{\log
\varepsilon }{a_{r}}\big).
\end{eqnarray*}Since\begin{eqnarray*}
\lim_{n\rightarrow \infty }\mathbf{P}\big(\tfrac{S_{\tau _{r}}-S_{r}}{a_{r}}
&\geq &\tfrac{\log \varepsilon }{a_{r}}\big) \\
&=&\mathbf{P}\big(\inf_{t\in \lbrack 0,1]}L_{t}-L_{1}\geq 0\big)=0,
\end{eqnarray*}%
where $L_{t}$ is\ a Levy process generated dy the stabe distribution with
index $\alpha $, we have\begin{equation*}
\lim_{n\rightarrow \infty }\mathbf{P}(e^{S_{\tau _{r}}-S_{r}}\geq
\varepsilon \mid \tau _{n-m}=n-m)=0.
\end{equation*}

Combaning this estimate with Lemmas \ref{le51} and \ref{le52} we conclude
that for all $\varepsilon >0$ there is a natural number $a$ such that for
sufficiently large $n$ \begin{equation}
\mathbf{P}\big(|e^{S_{\tau _{r}}-S_{r}}\tilde{Z}_{r}-\beta _{a,r}|\geq
\varepsilon /2\mid \tau _{n-m}=n-m\big)\leq \varepsilon /2.  \label{i1}
\end{equation}

Moreover, from Lemma \ref{leQ} we see that $\beta _{a,r},$ conditioned on
$\tau _{n-m}=n-m$, converges in distribution for every $a$. This implies that
$e^{S_{\tau _{r}}-S_{r}}\tilde{Z}_{r},$ conditioned on $\tau _{n-m}=n-m$
converges in distribution. According to Lemma \ref{le41} there is a $\delta
>0$ such that \begin{equation}
\mathbf{P}^{+}\Big(e^{-S_{a}}\sum_{1\leq i\leq a}\tilde{Z}_{a}^{i}<\delta
\Big)<\varepsilon /2\   \label{hvost}
\end{equation}if $a$ \ is sufficiently large.

For sufficiently small $\varepsilon >0$ let \begin{equation*}
\delta ^{\prime }=\delta -\varepsilon /2>0.
\end{equation*}In view of (\ref{i1}) and (\ref{hvost})\begin{eqnarray*}
\lim_{n\rightarrow \infty }\mathbf{P}(\beta _{a,r} &<&\delta ^{\prime }\mid
\tau _{n-m}=n-m) \\
&\leq &\lim_{n\rightarrow \infty }\mathbf{P}(e^{S_{\tau
_{r}}-S_{r}}\tilde{Z}_{r}<\delta ^{\prime }+|e^{S_{\tau _{r}}-S_{r}}\tilde{Z}_{r}-\beta
_{a,r}|\mid \tau _{n-m}=n-m) \\
&\leq &\lim_{n\rightarrow \infty }\mathbf{P}(e^{S_{\tau
_{r}}-S_{r}}\tilde{Z}_{r}<\delta \mid \tau _{n-m}=n-m)+\varepsilon /2\leq \varepsilon .
\end{eqnarray*}Thus\begin{equation*}
\mathbf{P}(\beta _{a,r}<\delta ^{\prime }\mid \tau _{n-m}=n-m)\leq
\varepsilon
\end{equation*}%
if $n$ is sufficiently large. Therefore the limiting distribution of
$e^{S_{\tau _{r}}-S_{r}}\tilde{Z}_{r}$ conditioned on $\tau _{n-m}=n-m$ has
no atom in zero. Applying Theorem \ref{trest} with $d_{n}=n-r_{n}$ , we
obtain the needed result. The theorem is proved.


\bigskip

\begin{thebibliography}{99}
\bibitem{agkv} \textsc{Afanasyev V.\thinspace I., Geiger J., Kersting G.,
Vatutin V.\thinspace A.} (2005). Criticality for branching processes in
random environment. \textsl{Ann. Probab.} \textbf{33}, 645--673.

\bibitem{ABKVINTER} \textsc{Afanasyev V.I., Boinghoff Ch., Kersting G.,
Vatutin V.A.} (2014). Conditional limit theorems for intermediately
subcritical branching processes in random environment. \textsl{Ann. Inst. H.
Poincarerob Probab. Statist.}, \textbf{50}:2, 602-627.

\bibitem{at} \textsc{Athreya K.B., Karlin, S.} (1971). On branching
processes with random environments: I, II. \textsl{Ann. Math. Stat.
}SUM{447}\textbf{42}, 1499--1520, 1843--1858.

\bibitem{be} \textsc{Bertoin J.} (1996). \textsl{L\'{e}vy processes}.
Cambridge University Press, Cambridge.

\bibitem{bedo} \textsc{Bertoin J., Doney, R.A.} (1994). On conditioning a
random walk to stay non-negative. \textsl{Ann. Probab. } \textbf{22}, 2152
-- 2167.

\bibitem{bi} \textsc{Bingham N.H., Goldie C.M., Teugels J.L.} (1987).
\textsl{Regular variation}. Cambridge University Press, Cambridge.

\bibitem{ca} \textsc{Caravenna F., Chaumont L}. (2008). Invariance
principles for random walks conditioned to stay positive. \textsl{Ann. I.H.
Poincar\'{e} (B). } \textbf{44}, 170--190.

\bibitem{chau} \textsc{Chaumont L.} (1997). Excursion normalis\'{e}e,
m\'{e}andre et pont pour les processus de L\'{e}vy stables. \textsl{Bull. Sci.
Math. } \textbf{121}, 377--403.

\bibitem{chau91} \textsc{Chauvin B., Rouault A., Wakolbinger A. } (1991).
Growing conditioned trees. \textsl{\ Stoch. Proc. Appl.} \textbf{39},
117--130.

\bibitem{ka77} \textsc{Kallenberg O.} (1977). Stability of critical cluster
fields.. \textsl{Math. Nachr.}~\textbf{77}, 7--43.

\bibitem{lpp} \textsc{Lyons R., Pemantle R., Peres Y.} (1995). Conceptual
proofs of $L\log L$ criteria for mean behavior of branching processes.
\textsl{Ann. Probab.} \textbf{23}, 1125--1138.

\bibitem{ne} \textsc{Neveu J.} (1986). Erasing a branching tree.
\textsl{Adv. Apl. Probab.} \textbf{suppl.}, 101--108.

\bibitem{sm} \textsc{Smith W.L., Wilkinson, W.E.} (1969). On branching
processes in random environments. \textsl{Ann. Math. Stat. } \textbf{40},
814--827.

\bibitem{dy04} \textsc{Vatutin V.A., Dyakonova E.E.} (2004). Galton-Watson
branching processes in random environment. I: limit theorems. \textsl{Theory
Probab. Appl.}, \textbf{48}, 314--336.

\bibitem{dy05} \textsc{Vatutin V.A., Dyakonova E.E.} (2005). Galton-Watson
branching processes in random environment. II: finite-dimensional
distributions. \textsl{Theory Probab. Appl.}, \textbf{49}, 275--308.

\bibitem{VaDy17} \textsc{Vatutin V., Dyakonova E. }(2017). Path to survival
for the critical branching processes in a random environment. \textsl{J.
Appl. Probab.}, \textbf{54}:2, 588-602.

\bibitem{VaDy19} \textsc{Vatutin V., Dyakonova E. }(2019).\ The initial
evolution stage of a weakly subcritical branching process in random
environment. \textsl{J. Appl. Probab.}, \textbf{64}:4, 535-552.

\bibitem{VaKe17} \textsc{Vatutin V., Kersting G. }(2017). Discrete time
branching processes in random environment, Wiley, New Jersey, USA, 306 pp.

\bibitem{Fel2} \textsc{Feller W.}\ (1967). An Introduction to Probability
Theory and its Applications. V. 2, Willey, New
York-London-Sydney-Toronto,1971.
\end{thebibliography}
\end{document}